\documentclass[12pt]{article}
\usepackage{amssymb, amsmath}

\newtheorem{theorem}{Theorem}[section]
\newtheorem{lemma}{Lemma}[section]
\newtheorem{definition}{Definition}[section]

\newcommand{\n}{\nonumber}

\newcommand{\si}{\sigma_R (|x|)}

\newcommand{\s}{\sigma}

\newcommand{\bb}{\begin{equation}}
\newcommand{\ee}{\end{equation}}
\newcommand{\bq}{\begin{eqnarray}}
\newcommand{\eq}{\end{eqnarray}}
\newcommand{\bqn}{\begin{eqnarray*}}
\newcommand{\eqn}{\end{eqnarray*}}

\begin{document}
\title{ On the nonexistence of  time dependent global weak solutions to
the compressible Navier-Stokes equations}
\author{ Dongho Chae\\
 Department of Mathematics\\
  Sungkyunkwan
University\\
 Suwon 440-746, Korea\\
 e-mail: {\it chae@skku.edu }}
 \date{}
\maketitle
\begin{abstract}
In this paper we prove the nonexistence of global weak solutions to
the compressible Navier-Stokes equations for the isentropic gas in
$\Bbb R^N, N\geq 3,$ where the pressure law given by
$p(\rho)=a\rho^{\gamma}, $ $a>0, 1<\gamma \leq \frac{N}{4}+\frac12$.
In this case if the initial data satisfies $\int_{\Bbb R^N} \rho_0
(x)v_0 (x)\cdot x\, dx >0$, then there exists no  finite energy
global weak solution which satisfies  the integrability conditions $
\rho |x|^2 \in L^1_{\mathrm{loc}} (0, \infty; L^1 (\Bbb R^N))$ and $
v\in
L^1_{\mathrm{loc}} (0, \infty; L^{\frac{N}{N-1}} (\Bbb R^N))$.   
{\bf keywords:} compressible Navier-Stokes equations, nonexistence
of the global weak solutions
\end{abstract}
\section{Introduction and Main Theorems}
\setcounter{equation}{0} We are concerned  on the compressible
Navier-Stokes equations on $\Bbb R^N$, $N\geq 1$.
  $$
(NS)\left\{\aligned &\partial_t \rho + \mbox{div}(\rho v) = 0, \\
& \partial_t (\rho v) + \mbox{div}(\rho v \otimes v) = -\nabla
p+\mu \Delta v +(\mu +\lambda ) \nabla \mathrm{div }\, v+f, \\
& p(\rho)=a \rho ^\gamma, \,\gamma >1, a>0
\endaligned \right.
$$
The system (NS) describes  isentropic viscous compressible gas
flows, and $\rho, v$, $p$ and $f$ denote the density, velocity,
 pressure and the external force respectively.  We treat  the viscous case $\mu
>0, \lambda>0$.  For the surveys of the known mathematical theories of
the system (NS) we refer to \cite{lio, fei, nov}. In this paper we
are concerned on the case where $1<\gamma \leq N/4+1/2,$ and $N\geq
3$. We note that for this range of $\gamma, N$ there exists no known
existence results for the global weak solutions(see \cite{lio,fei}
for previous global existence results of the weak solutions). In
this case of our aim  is to prove nonexistence of global weak
solutions to (NS) satisfying suitable integrability conditions for
the solutions and for some class of initial data. This result can be
regarded as a time dependent generalization of the corresponding
stationary results obtained previously in \cite{cha}. The finite
energy weak solution of (NS) is defined as follows.
  \begin{definition}
  Let $(\rho_0,  v_0)$ satisfies
 \bb
  0\leq\rho_0 \in L^1_{\mathrm{loc}} (\Bbb R^N), \quad \rho_0 v_0 \in
  L^1_{\mathrm{loc}}(\Bbb R^N),\quad \rho_0 ^{\gamma} +\rho_0 |v_0 |^2 \in L^1 (\Bbb R^N).
  \ee
  We say that the pair $(\rho, v)$  is a finite energy weak solution of the system
  (NS) with the initial data $(\rho, v )$ if it satisfies
\bb\label{AA}
  0\leq\rho  \in L^1_{\mathrm{loc}} (\Bbb R^N\times (0, \infty)), \quad \rho v \in
  L^1_{\mathrm{loc}}(\Bbb R^N\times(0, \infty)),
\ee
  \bb\label{BB}
    \rho ^{\gamma} +\rho |v |^2
  \in L^\infty_{\mathrm{loc}}(0, \infty;L^1 (\Bbb R^N)), \quad v\in L^2 _{loc} (0,\infty; \dot{H} (\Bbb
  R^N)),
  \ee
  and satisfies (NS) in the sense of distribution, namely
  \bq\label{1}
 && \xi(0)\int_{\Bbb R^N} \rho_0 (x) \psi (x)dx+
 \int_0 ^\infty \int_{\Bbb R^N} \rho (x,t) \psi (x) \xi ' (t)dxdt\n\\
 &&\qquad+\int_0 ^\infty\int_{\Bbb R^N} \rho v(x,t)\cdot \nabla \psi (x)\xi(t)\,dx=0 \qquad \forall \psi
 \in C_0 ^\infty (\Bbb R^N),\xi \in C_0 ^1([0, \infty)),\n\\
 \ \\
 \label{2}
 &&\xi(0)\int_{\Bbb R^N} \rho_0 (x)v_0(x)\cdot \phi (x)dx+
 \int_0 ^\infty \int_{\Bbb R^N} \rho (x,t)v(x,t)\cdot \phi (x) \xi ' (t)dxdt\n \\
 &&\qquad+\int_0 ^\infty\int_{\Bbb R^N} \rho(x,t) v(x,t)\otimes v(x,t) :\nabla \phi (x)\xi(t)\,dxdt\n \\
 &&\qquad=
 -\int_0 ^\infty\int_{\Bbb R^N} p(x,t)\, \,\mathrm{div }\,
 \phi(x)\xi(t)
 \,dxdt-\mu\int_0 ^\infty\int_{\Bbb R^N} v(x,t)\cdot \Delta \phi(x)\xi(t)\, dxdt\n \\
 &&\qquad -(\mu +\lambda
 )\int_0 ^\infty\int_{\Bbb R^N} v(x,t)\cdot \nabla \mathrm{div}\, \phi (x)\xi(t)\,
 dxdt-\int_0 ^\infty\int_{\Bbb R^N} f\cdot \phi(x)\xi(t)\, dxdt \n \\
 &&\qquad\qquad\qquad
\quad \forall \phi
 \in [C_0 ^\infty (\Bbb R^N)]^N,\xi \in C_0 ^1([0, \infty)).
 \eq
\bb
 \label{4}
 p(\rho)=a\rho^\gamma \quad \mbox{almost everywhere in $\Bbb R^N\times (0, \infty)$},
  \ee
   Finally, we impose the energy inequality of the following form:
\bq\label{energy} &&\frac{d}{dt} \mathcal{E}(t)+ \int_{\Bbb R^N}
(\mu |\nabla v|^2 +(\mu+\lambda) | \mathrm{div}\, v|^2 )
dx\leq  0 \quad \mbox{for almost every $ t\geq 0$}, \n \\
&&\quad\mbox{where} \qquad \mathcal{E}(t):=\int_{\Bbb R^N}
\left[\frac12 \rho |v|^2 +\frac{a \rho^{\gamma} }{\gamma-1}\right]
dx.
 \eq
  \end{definition}
  In the above definition we followed closely \cite{nov}. The
  following is our main theorem.

\begin{theorem}
Let $N\geq 3$, $1<\gamma \leq \frac{N}{4}+\frac12$,  $\mu>0 ,
\mu+\lambda
>0$, and the external force $f\in [L^1_{loc}(\Bbb R^N\times[0, \infty))]^N$
with $\mathrm{div}\, f=0$ is given.  Let the initial
 data $(\rho_0, v_0)$  satisfy
 \bb\label{7}
  \mathcal{E}(0) <\infty,\quad\int_{\Bbb R^N} \rho_0 (x)|v_0 (x)| |x|\,dx <\infty,
 \ee
  \bb\label{11b}
 \int_{\Bbb R^N} \rho_0 (x)v_0 (x)\cdot x\, dx >0.
 \ee
Then, there exists no finite energy  global weak solution to (NS)
such that
  \bb\label{30a}
 \rho |x|^2\in L^1_{\mathrm{loc}}(0, \infty; L^1(\Bbb R^N)),\quad
 v\in L^1_{\mathrm{loc}} (0,
\infty; L^{\frac{N}{N-1}} (\Bbb R^N)).
 \ee
\end{theorem}
{\em Remark 1.1 } An immediate consequence of the above theorem is
that for initial data $(\rho_0, v_0 )$ satisfying (\ref{7}) and
(\ref{11b})  there exists $T_* <\infty$ such that
$$ \int_0 ^{T_*}\int_{\Bbb R^N}\left[\rho(x,t) |x|^2 + |v(x,t)|^{\frac{N}{N-1}} \right]dx dt
=\infty,
$$
where $(\rho(x,t), v(x,t))$ is a local in time (classical or weak)
solution on the time interval $[0, T_*)$, if (NS) is at least
locally well-posed in the functional setting of Definition 1.1 for
the case considered above.
\section{Proof of Theorem 1.1}
\setcounter{equation}{0}

 In order to prove Theorem 1.1 we shall use the following lemma, which is proved in \cite{guo}.
\begin{lemma}
Suppose $(\rho, v)$ is a finite energy global weak solution to
$(NS)$ with the setting given by Theorem 1.1, then
 \bb\label{guo}
 \int_0 ^\infty \int_{\Bbb R^N}
 \frac{ \rho(x,t)(1+|x|^2)^{\frac{N+2}{4\gamma}}}{t^2} dxdt \leq C \mathcal{E}(0).
 \ee
\end{lemma}
Since $\frac{N+2}{4\gamma}\geq 1$ in our setting of Theorem 1.2, one
immediate consequence of (\ref{guo}) is the following fact
  \bb\label{key}
\lim_{\tau\to \infty} \int_{\tau}^{2\tau}\int_{\Bbb R^N} \frac{\rho
(x,t)|x|^2}{1+t^2}\,dxdt=0.
 \ee
 Indeed, using (\ref{guo}), we deduce
 $$
\lim_{\tau\to \infty} \int_{\tau}^{2\tau}\int_{\Bbb R^N} \frac{\rho
(x,t)|x|^2}{1+t^2} \,dxdt\leq\lim_{\tau\to \infty}
\int_{\tau}^{2\tau}\int_{\Bbb R^N}
 \frac{ \rho(x,t)(1+|x|^2)^{\frac{N+2}{4\gamma}}}{t^2} dxdt=0,
 $$
 where the last step follows from the dominated convergence
 theorem.\\
 \ \\
 \noindent{\bf Proof of Theorem 1.2 }
Suppose there exists a global weak solution $(\rho, v)$ satisfying
(\ref{1})-(\ref{energy}).
 Let us
consider a radial cut-off function $\sigma\in C_0 ^\infty(\Bbb R^N)$
such that
 \bb\label{z12}
   \sigma(|x|)=\left\{ \aligned
                  &1 \quad\mbox{if $|x|<1$}\\
                     &0 \quad\mbox{if $|x|>2$},
                      \endaligned \right.
 \ee
and $0\leq \sigma  (x)\leq 1$ for $1<|x|<2$.  For $R>0$ we define
 \bb\label{z33}
\varphi_R (x)=\frac12|x|^2 \s
\left(\frac{|x|}{R}\right)=\frac12|x|^2\s_R (|x|)\in C_0 ^\infty
(\Bbb R^N).
 \ee
 We also introduce temporal cut-off function $\eta \in C^\infty_0 ([0, \infty))$ as follows.
\bb\label{z34}
   \eta(t)=\left\{ \aligned
                  &1 \quad\mbox{if $0\leq t<1$}\\
                     &0 \quad\mbox{if $t>2$},
                      \endaligned \right.
 \ee
and for all $\tau >0$ we set
 \bb\label{z35}
 \eta_\tau(t)=\eta\left(\frac{t}{\tau}\right).
 \ee
Substituting $\phi(x)=\nabla \varphi_R (x), \xi(t)=\eta_\tau(t)$
into (\ref{2}), we obtain
 \bq\label{z36}
 \lefteqn{0=\int_{\Bbb R^N} \rho_0 (x)v_0(x)\cdot x\si dx+\frac{1}{2 R}\int_{\Bbb R^N}
  \rho_0 (x)v_0(x)\cdot x |x|\s'\left(\frac{|x|}{R}\right)
 dx}\hspace{.0in}\n \\
 &&+\int_0 ^\infty \int_{\Bbb R^N} \rho (x,t)v(x,t) \cdot
 \nabla \varphi_R (x)
 \eta_\tau ^{\prime}(t) dxdt \n \\
&&+\int_0 ^\infty\int_{\Bbb R^N} \rho (x,t) |v(x,t)|^2 \si \eta_\tau(t) \,dx dt\n\\
&&+\frac{1}{2R}\int_0 ^\infty\int_{\Bbb R^N} \rho(x,t)  \s'
\left(\frac{|x|}{R}\right) \frac{(v(x,t)\cdot x)^2}{|x|}\eta_\tau(t)
\,dxdt \n \\
 &&+ \frac{1}{2R}\int_0 ^\infty\int_{\Bbb R^N}
\rho(x,t)|v(x,t)|^2 |x| \s'\left(\frac{|x|}{R}\right)\eta_\tau(t)
\,dxdt \n \\
  &&+\frac{1}{2R^2}\int_0 ^\infty\int_{\Bbb R^N} \rho(x,t)(v(x,t)\cdot x)^2}{
\s^{\prime\prime} \left(\frac{|x|}{R}\right) \eta_\tau(t) \,dx dt\n \\
 &&+ N\int_0 ^\infty \int_{\Bbb R^N}p(x,t)\sigma_R (|x|)\eta_\tau(t)
  \, dxdt\n \\
&&+ \frac{2}{R}\int_0 ^\infty\int_{\Bbb R^N}p(x,t) |x|
 \s' \left(\frac{|x|}{R}\right)\eta_\tau(t)\, dxdt\n \\
&&+ \frac{N-1}{2R}\int_0 ^\infty\int_{\Bbb R^N}p(x,t)|x|\s'
\left(\frac{|x|}{R}\right) \eta_\tau(t) \, dxdt\n \\
&&+ \frac{1}{2R^2}\int_0 ^\infty\int_{\Bbb R^N}p(x,t)|x|^2
\s^{\prime\prime} \left(\frac{|x|}{R}\right)\eta_\tau(t) \,
dxdt \n \\
&&+(2\mu +\lambda)\int_0 ^\infty\int_{\Bbb R^N} v\cdot\nabla \Delta
(|x|^2
 \s \left(\frac{|x|}{R}\right)\eta_\tau(t)\, dxdt,\n \\
  &&:=I_1+\cdots +I_{12}.
 \eq
On the other hand, substituting $\phi(x)= \nabla\varphi_R (x)$, $\xi
(t)=\eta_\tau'(t)$
 into (\ref{1}),  we find that (note that $\xi (0)=\eta_\tau'(0)=0$)
 \bq\label{z37}
 I_3&=&\int_0 ^\infty\int_{\Bbb R^N} \rho v(x,t)\cdot
\nabla
\varphi_R (x)\eta_{\tau} '(t)\,dxdt\n \\
&=&- \int_0 ^\infty \int_{\Bbb R^N} \rho (x,t) \varphi_R (x)
 \eta^{\prime\prime}_\tau
 (t)dxdt\n \\
 &=&-\frac12 \int_0 ^\infty \int_{\Bbb R^N} \rho (x,t)|x|^2\sigma_R (|x|)\eta^{\prime\prime}_\tau
 (t) dxdt\n \\
 &&\to-\frac12\int_0 ^\infty \int_{\Bbb R^N} \rho (x,t)
|x|^2\eta^{\prime\prime}_\tau
 (t) dxdt
 \eq
 as $R\to \infty$ by the dominated convergence theorem.
  We  also have
  \bb\label{z38}
  I_4 \to \int_0 ^\infty\int_{\Bbb R^N} \rho(x,t)|v(x,t)|^2  \eta_\tau (t)\,dxdt
 \ee
  as $R\to \infty$.
Similarly,
 \bb\label{z38a}
 I_1\to \int_{\Bbb R^N} \rho_0 (x)v_0(x)\cdot x dx,
 \ee
 and
  \bb\label{z39}
  I_8\to N\int_0 ^\infty \int_{\Bbb R^N}p(x,t)\eta_\tau (t)
 \, dxdt
\ee
as $R\to \infty$.
 For $I_5, I_6$ we estimate
 \bq\label{z40}
  |I_5 |+|I_6|&\leq& \int_0 ^{2\tau}\int_{R<|x|<2R} \rho(x,t)|v(x,t)|^2\left|\s'
\left(\frac{|x|}{R}\right)\right|
  \frac{|x|}{R}dxdt\n \\
  &\leq &2 \sup_{1<s<2} |\s'(s)|
\int_0 ^{2\tau}\int_{R<|x|<2R}\rho(x) |v(x,t)|^2
 \, dxdt
\to 0\n \\
  \eq
 as $R\to \infty$.
Similarly
 \bb\label{z40a}
  |I_2|\leq \int_{R<|x|<2R} \rho_0 (x) |x| dx \to 0,
\ee
 and
 \bq\label{z41}
  |I_7|&\leq&\frac12\int_0 ^{2\tau}\int_{R<|x|<2R}\frac{|x|^2}{R^2}
   \rho(x,t)|v(x,t)|^2\left|\s^{\prime\prime}
  \left(\frac{|x|}{R}\right)\right| \,dxdt\n \\
  &\leq&2\sup_{1<s<2} |\s^{\prime\prime}(s)|
 \int_0 ^{2\tau}\int_{R<|x|<2R}\rho(x,t)|v(x,t)|^2
 \,dx \to 0\n \\
  \eq
   as $R\to \infty$. The estimates for $I_9,I_{10}$ and $I_{11}$ are
   similar to the above, and we find
   \bq\label{z42}
   |I_9|&\leq &2 \int_0 ^{2\tau} \int_{R<|x|<2R}|p(x,t)| \frac{|x|}{R}
 \left|\s' \left(\frac{|x|}{R}\right)\right|  \, dxdt\n \\
 &\leq& 4\sup_{1<s<2} |\s' (s)|
\int_0 ^{2\tau}\int_{R<|x|<2R}|p(x,t)|\,dxdt \to 0,\n \\
  \eq
\bq\label{z43}
   |I_{10}|&\leq & \frac{N-1}{2R}\int_0 ^{2\tau}\int_{R<|x|<2R}|p(x,t)||x|\left|\s'
\left(\frac{|x|}{R}\right)\right|  \, dxdt\n \\
 &\leq& (N-1)\sup_{1<s<2} |\s' (s)|
\int_0 ^{2\tau}\int_{R<|x|<2R}|p(x,t)|dxdt\to 0,\n \\
  \eq
  and
  \bq\label{z44}
 |I_{11}|&\leq&\frac{1}{2R^2} \int_0 ^{2\tau}\int_{\Bbb
 R^N}|p(x,t)||x|^2
\left|\s^{\prime\prime}\left(\frac{|x|}{R}\right)\right|
 \, dxdt \n\\
  &\leq& 2\sup_{1<s<2}
|\s^{\prime\prime} (s)|
 \int_0 ^{2\tau}\int_{R<|x|<2R}|p(x,t)|\,dxdt \to 0\n \\
  \eq
  as $R\to \infty$ respectively. Now we show the vanishing of the
  viscosity term as $R\to \infty$.
  This follows from the estimates,
 \bq\label{z45}
  |I_{12}|&=&(2\mu +\lambda ) \left|\int_0 ^\infty\int_{\Bbb R^N} v\cdot\nabla \Delta (|x|^2
 \s \left(\frac{|x|}{R}\right)\eta_\tau(t)\,dxdt\right|\n \\
&\leq& (2\mu +\lambda) \left|\int_0 ^\infty\int_{\Bbb R^N}
(N+5)\left[\frac{(v\cdot x)}{R|x|}\s ' \left(\frac{|x|}{R}\right)+
\frac{(v\cdot x)}{R^2}\s^{\prime\prime}
\left(\frac{|x|}{R}\right)\right]\eta_\tau(t)\, dxdt\right|\n \\
&&\qquad+(2\mu +\lambda)\left|\int_0 ^\infty\int_{\Bbb R^N}
\frac{|x|(v\cdot x)}{R^3}\s^{\prime\prime\prime}
\left(\frac{|x|}{R}\right)\eta_\tau(t)\,dx dt\right|\n \\
&\leq&  \frac{C}{R} \int_0 ^{2\tau}\int_{R\leq |x|\leq 2R} |v(x,t)|
\, dx dt\n \\
&\leq& C \int_0 ^{2\tau}\left(\int_{R\leq  |x|\leq 2R}
|v(x,t)|^{\frac{N}{N-1}} dx\right)^{\frac{N-1}{N}} dt\to 0
 \eq
   as $R\to \infty$.
   In summary,  passing  $R\to \infty$ in (\ref{z36}),
  we obtain
\bq\label{z46}
 &&\frac12\int_0 ^\infty \int_{\Bbb R^N} \rho (x,t)
|x|^2\eta^{\prime\prime}_\tau
 (t) dxdt=\frac12\int_0 ^\infty\int_{\Bbb R^N} \rho(x,t) |v(x,t)|^2\eta_\tau (t)
 \,dxdt \n \\
&&\qquad+N\int_0 ^\infty\int_{\Bbb R^N}p(x,t)\eta_\tau
(t)dxdt+\int_{\Bbb R^N} \rho_0 (x)v_0(x)\cdot x\,
 dx,
 \eq
 and therefore
\bb\label{z47}
 \int_{\Bbb R^N} \rho_0 (x)v_0(x)\cdot x\,dx \leq \frac12\int_0 ^\infty \int_{\Bbb R^N} \rho (x,t)
|x|^2\eta^{\prime\prime}_\tau
 (t) dxdt
 \ee
for any $\tau >0$. By (\ref{key}) we have
  \bq\label{z48}
   \lefteqn{\left|\int_0 ^\infty \int_{\Bbb R^N} \rho (x,t)
|x|^2\eta^{\prime\prime}_\tau
 (t) dxdt\right|\leq \frac{1}{\tau^2} \int_{\tau}^{2\tau}
 \int_{\Bbb R^N} \rho (x,t)
|x|^2\left|\eta^{\prime\prime} \left(\frac{t}{\tau}\right)\right|
  dxdt }\hspace{.4in}\n \\
  &&\leq\frac{1+4\tau^2}{\tau^2} \sup_{1<s<2} |\eta^{\prime\prime}(s)|
  \int_{\tau}^{2\tau}\int_{\Bbb R^N} \frac{\rho (x,t)|x|^2}{1+t^2}
dxdt\to 0
 \eq
  as $\tau\to \infty$.
 Combining (\ref{z48}) with (\ref{z47}), we obtain the following necessary condition for the global existence of
 weak solution,
 \bb \int_{\Bbb R^N} \rho_0 (x)v_0(x)\cdot
x\, dx\leq 0.
 \ee
This proves the theorem. $\square$\\

$$\mbox{ \bf Acknowledgements} $$
 This
work was supported partially by  KRF Grant(MOEHRD, Basic Research
Promotion Fund).

\end{document}